\documentclass[a4paper,12pt]{amsart}
\usepackage{amsfonts}
\usepackage{amssymb}
\usepackage{ifthen}
\usepackage{graphicx}
\usepackage{enumitem}
\usepackage{amsmath}
\nonstopmode \numberwithin{equation}{section}
\usepackage[usenames]{color}
\newtheorem{thm}{Theorem}[section]
\newtheorem{lem}{Lemma}[section]
\newtheorem{cor}[thm]{Corollary}
\newtheorem{prop}[thm]{Proposition}

\theoremstyle{definition}
\newtheorem{mlem}{Main lemma}[section]
\newtheorem{assertion}{Assertion}[section]
\newtheorem{cl}{Claim}[section]
\newtheorem{ca}{Case}[section]
\newtheorem{sca}{Subcase}[section]
\newtheorem{scl}{Subclaim}[section]
\newtheorem{conj}[thm]{Conjecture}
\newtheorem{fact}{Fact}[section]
\newtheorem{defn}{Definition}[section]
\newtheorem{prob}{Problem}
\newtheorem{ques}[thm]{Question}
\newtheorem{rem}{Remark}[section]
\newtheorem{exam}{Example}[section]

\numberwithin{equation}{section}
\DeclareMathOperator*{\esssup}{ess\,sup}
\DeclareMathOperator*{\essinf}{ess\,inf}
\DeclareMathOperator*{\dlimsup}{lim\,sup}

\newcounter {own}
\def\theown {\thesection       .\arabic{own}}

\newlist{steps}{enumerate}{1}
\setlist[steps,1]{
  leftmargin=*,
  label=\textbf{Step \arabic*}.,
  ref=Step~\arabic*,
}

\newenvironment{pf}[1][]{%
 \vskip 3mm
 \noindent
 \ifthenelse{\equal{#1}{}}%
  {{\slshape Proof. }}%
  {{\slshape #1.} }%
 }%
{\qed\bigskip}

\newcounter{alphabet}
\newcounter{tmp}

\newcounter{alphabet2}

\makeatletter
\newcommand{\Ref}[1]{\@ifundefined{r@#1}{}{\setcounter{tmp}{\ref{#1}}\Alph{tmp}}}
\makeatother

\newcommand{\IR}{{\mathbb R}}

\newcommand{\ID}{{\mathbb D}}

\def\be{\begin{equation}}
\def\ee{\end{equation}}

\newcommand{\ben}{\begin{enumerate}}
\newcommand{\een}{\end{enumerate}}

\newcommand{\blem}{\begin{lem}}
\newcommand{\elem}{\end{lem}}
\newcommand{\bthm}{\begin{thm}}
\newcommand{\ethm}{\end{thm}}
\newcommand{\bcor}{\begin{cor}}
\newcommand{\ecor}{\end{cor}}
\newcommand{\beg}{\begin{exam}}
\newcommand{\eeg}{\end{exam}}
\newcommand{\begs}{\begin{examples}}
\newcommand{\eegs}{\end{examples}}
\newcommand{\bdefe}{\begin{defn}}
\newcommand{\edefe}{\end{defn}}
\newcommand{\bprob}{\begin{prob}}
\newcommand{\eprob}{\end{prob}}
\newcommand{\bques}{\begin{ques}}
\newcommand{\eques}{\end{ques}}
\newcommand{\bei}{\begin{itemize}}
\newcommand{\eei}{\end{itemize}}
\newcommand{\bcon}{\begin{conj}}
\newcommand{\econ}{\end{conj}}
\newcommand{\bop}{\begin{op}}
\newcommand{\eop}{\end{op}}

\newcommand{\bas}{\begin{assertion}}
\newcommand{\eas}{\end{assertion}}

\newcommand{\bfa}{\begin{fact}}
\newcommand{\efa}{\end{fact}}

\newcommand{\bca}{\begin{ca}}
\newcommand{\eca}{\end{ca}}
\newcommand{\bsca}{\begin{sca}}
\newcommand{\esca}{\end{sca}}

\newcommand{\bcl}{\begin{cl}}
\newcommand{\ecl}{\end{cl}}

\newcommand{\bmlem}{\begin{mlem}}
\newcommand{\emlem}{\end{mlem}}

\newcommand{\bscl}{\begin{scl}}
\newcommand{\escl}{\end{scl}}

\newcommand{\bcons}{\begin{conjs}}
\newcommand{\econs}{\end{conjs}}

\newcommand{\bprop}{\begin{prop}}
\newcommand{\eprop}{\end{prop}}

\newcommand{\br}{\begin{rem}}
\newcommand{\er}{\end{rem}}
\newcommand{\brs}{\begin{rems}}
\newcommand{\ers}{\end{rems}}
\newcommand{\bo}{\begin{obser}}
\newcommand{\eo}{\end{obser}}
\newcommand{\bos}{\begin{obsers}}
\newcommand{\eos}{\end{obsers}}
\newcommand{\bpf}{\begin{pf}}
\newcommand{\epf}{\end{pf}}
\newcommand{\ba}{\begin{array}}
\newcommand{\ea}{\end{array}}
\newcommand{\beq}{\begin{eqnarray}}
\newcommand{\beqq}{\begin{eqnarray*}}
\newcommand{\eeq}{\end{eqnarray}}
\newcommand{\eeqq}{\end{eqnarray*}}

\newcounter{minutes}
\divide\time by 60
\newcounter{hours}
\multiply\time by 60 \addtocounter{minutes}{-\time}

\begin{document}

\bibliographystyle{amsplain}

\title{Norm estimates of the partial derivatives for harmonic mappings and harmonic quasiregular mappings}


\author{Jian-Feng Zhu}
\address{Jian-Feng Zhu,  Department of Mathematics,
Shantou University, Shantou, Guangdong 515063, People's Republic of China and School of Mathematical Sciences,
Huaqiao University,
Quanzhou 362021, People's Republic of China.
} \email{flandy@hqu.edu.cn}

\date{}
\subjclass[2000]{Primary 30C55, 30C62}
\keywords{Hardy norm, Bergman norm, harmonic mapping, quasiregular mapping, Poisson kernel.\\
}

\begin{abstract}
Suppose $p\geq1$, $w=P[F]$ is a harmonic mapping of the unit disk $\ID$ satisfying $F$ is
absolutely continuous and $\dot{F}\in L^p(0, 2\pi)$, where $\dot{F}(e^{it})=\frac{\mathrm{d}}{\mathrm{d}t}F(e^{it})$.
In this paper, we obtain Bergman norm estimates of the partial derivatives for $w$, i.e.,
$\|w_z\|_{L^p}$ and $\|\overline{w_{\bar{z}}}\|_{L^p}$, where $1\leq p<2$.
Furthermore, if $w$ is a harmonic quasiregular mapping of $\mathbb{D}$, then
we show that $w_z$ and $\overline{w_{\bar{z}}}$ are in the Hardy space $H^p$, where $1\leq p\leq\infty$.
The corresponding Hardy norm estimates, $\|w_z\|_{p}$ and $\|\overline{w_{\bar{z}}}\|_{p}$, are also obtained.
\end{abstract}
\thanks{}

\maketitle \pagestyle{myheadings} \markboth{Jian-Feng Zhu}{Norm estimates for harmonic and harmonic qr. mappings}


\section{Introduction}\label{sec-1}
In this paper, we mainly deal with planar harmonic mappings and planar quasiregular mappings.
For the convenient of stating our motivations and results,
we introduce the definitions of the Bergman norm, the Hardy norm and quasiregular mappings in $n$-dimensional.

Throughout this paper, we let $B(x, r)$ be the open ball in $\mathbb{R}^n$\, ($n\geq2$) with the radius $r$ and centered at $x$, denote by $\mathbb{B}^n$ the unit ball of $\mathbb{R}^n$, i.e., $\mathbb{B}^n=B(0, 1)$. Given $x\in \mathbb{B}^n$, we write $B_x=B(x, (1-|x|)/2)$. The boundary of $B(x, r)$ is denoted by $\mathbb{S}^{n-1}(x, r)$ and we write $\mathbb{S}^{n-1}=\mathbb{S}^{n-1}(0, 1)$.
For $n=2$, we let $\ID$ be the unit disk in the complex plane $\mathbb{C}$, and $\mathbb{T}$ the unit circle.
\subsection*{Bergman norm}
Denote by $L^p(\mathbb{B}^n)\ \,(1\leq p\leq\infty)$ the space of
measurable functions on $\mathbb{B}^n$ with finite integral
$$\|f\|_{L^p}=\left(\int_{\mathbb{B}^n}|f(x)|^p\mathrm{d}m(x)\right)^{\frac{1}{p}},\ \ \ 1\leq p<\infty,$$
where
$\mathrm{d}m(x)$ is the normalized Lebesgue measure on $\mathbb{B}^n$, i.e., $\int_{\mathbb{B}^n}\mathrm{d}m(x)=1$.
For the case $p=\infty$, we let $L^{\infty}(\mathbb{B}^n)$ denote the space of (essentially) bounded
functions on $\mathbb{B}^n$. For $f\in L^{\infty}(\mathbb{B}^n)$, we define
$$\|f\|_\infty=\esssup\{|f(x)|: x\in\mathbb{B}^n\}.$$

If in particular $n=2$, then we use $\mathrm{d}A(z)$ instead of $\mathrm{d}m(x)$ for the normalized Lebesgue measure, i.e., for $z=(x, y)\in \mathbb{R}^2$ or
$z=x+iy=re^{i\theta}\in\ID$, we write
$\mathrm{d}A(z) =\frac{1}{\pi}\mathrm{d}x\mathrm{d}y=\frac{1}{\pi}r\mathrm{d}r\mathrm{d}\theta$\,  (cf. \cite[Page 1]{Zhukh}).
The norm $\|f\|_{L^p}$ is called the {\it Bergman norm} of $f$\, (cf. \cite{Yamashita}) and
the space $L^{\infty}(\ID)$ is a Banach space with the above norm (cf. \cite[Page 2]{Zhukh}).

\subsection*{Hardy norm}
Let $f$ be an {\it analytic} function of $\ID$.
Following the notation of \cite{Du-hp}, the {\it integral means} of $f$ are defined as follows:
$$M_p(r, f)=\left\{\frac{1}{2\pi}\int_0^{2\pi}|f(re^{i\theta})|^p\mathrm{d}\theta\right\}^{1/p}, \ \ \ 0<p<\infty;$$
and
$$M_{\infty}(r, f)=\max\limits_{0\leq \theta\leq 2\pi}|f(re^{i\theta})|.$$
A function $f$ analytic in $\ID$ is said to be of class $H^p$\, $(0<p\leq\infty)$, if $M_p(r, f)$ remains bounded as $r$
tends to 1.

The norm
$$\|f\|_{p}=\lim\limits_{r\rightarrow1^-}M_p(r, f)$$
is called the {\it Hardy norm} of $f$, where $0<p\leq\infty$ (cf. \cite{Du-hp} and \cite{Yamashita}).

It is convenient also to define the analogous classes of {\it harmonic mappings}. A mapping $w(z)$ harmonic in $\ID$ is
said to be of class $h^p$\, $(0<p\leq \infty)$ if $M_p(r, w)$ is bounded. It is evident that $H^q\subset H^p$, if $0<p<q\leq\infty$, and
likewise for the $h^p$ spaces. Also, it is evident that $H^p\subseteq L^p(\ID)$ and $h^p\subseteq L^p(\ID)$, for all $p\geq1$.

Adopting the above classical definition, we say that a {\it quasiconformal} mapping (see the definition below) $f$ on $\mathbb{B}^n$\,  ($n\geq2$) belongs to the class $H^p$ provided (cf. \cite[Page 23]{Astala-Hp}) that
$$\|f\|_{p}=\sup\limits_{0<r<1}\left(\int_{\mathbb{S}^{n-1}}|f(r\omega)|^p\mathrm{d}\sigma(\omega)\right)^{1/p}<\infty,$$
where $\omega\in \mathbb{S}^{n-1}$ and $\mathrm{d}\sigma(\omega)$ is the normalized Lebesgue measure on $\mathbb{S}^{n-1}$.
According to Beurling's theorem, for a given quasiconformal mapping $f$, the radial limit
$$F(\omega)=\lim_{r\to 1^-}f(r\omega)$$
exists for a.e. $\omega\in \mathbb{S}^{n-1}$.
Define $\mathcal{M}(r, f):=\sup_{\omega\in\mathbb{S}^{n-1}}|f(r\omega)|$ for $0<r<1$. Then the {\it weighted Hardy space}, for $-1<\alpha<\infty$ and $0<p<\infty$, is defined as the class of all univalent functions for which (cf.  \cite[Page 1]{Koskela-19})
$$\int_0^1\mathcal{M}(r, f)^p(1-r)^{\alpha}\mathrm{d}r<\infty.$$
\subsection*{Poisson integral}
Suppose $w(z)=u(z)+iv(z)$\  ($z=x+iy$) is a complex-valued harmonic mapping of $\ID$.
Then, there exists analytic functions $g$ and $h$ defined on $\ID$  such that $w$ has the canonical
representation $w=h+\overline{g}$.
Also, every bounded harmonic mapping
$w$ defined on $\ID$ has the following representation
\be\label{Z1-eq1}
w(z)=P[F](z)=\int\limits_0^{2\pi} P_r(t-\theta)F(e^{it})\,\mathrm{d}t, \quad z=re^{i\theta}\in \ID,
\ee
where $F$ is a bounded integrable function defined on the unit circle $\mathbb{T}$,
and
$$P_r(t-\theta)=\frac{1}{2\pi}\frac{1-r^2}{1-2r\cos(t-\theta)+r^2},$$
denotes the {\it Poisson kernel}. We refer to \cite{Du-04} for more details and discussions on harmonic mappings.

For $F\in L^p(0, 2\pi)$, let
$$\|F\|_{L^p}=\left(\frac{1}{2\pi}\int_0^{2\pi}|F(e^{it})|^p\mathrm{d}t\right)^{1/p},\ \ \ \ 1\leq p<\infty.$$
If $p=\infty$, then we write
$$\|F\|_{\infty}:=\esssup\{|F(e^{it})|: t\in[0, 2\pi]\}.$$
It is known that if $w=P[F]$ is the {\it Poisson integral} of a function $F\in L^p(0, 2\pi)$, $1\leq p\leq\infty$, then
$w\in h^p$ and $M_p(r, w)\leq \|F\|_{L^p}$\,  (cf. \cite[Page 11]{Du-hp}).

\subsection*{Directional derivative and Jacobian}
The {\it formal derivatives} of a complex-valued function $w$ are defined by:
\beqq\label{eq1.1}
w_z=\frac{1}{2}\left(w_x-iw_y\right)\;\;\mbox{and}\;\;w_{\bar{z}}=\frac{1}{2}\left(w_x+iw_y\right),
\eeqq where $z=x+iy\in \ID$, and $x$, $y\in \IR$.

Assume that $z=re^{i\theta}\in\ID$, then the polar derivatives of $w$ are given as follows:
\be\label{5.3}
w_{\theta}(z)=i\big(zw_z(z)-\bar{z}w_{\bar{z}}(z)\big)\ \ \ \mbox{and}\ \ \ rw_r(z)=zw_z(z)+{\bar{z}}w_{\bar{z}}(z).
\ee
These show that $w_{\theta}(z)$ and $rw_r(z)$ are harmonic in $\ID$ and
\be\label{Oct-20-1} w_z(z)=\frac{e^{-i\theta}}{2}\left(w_r(z)-\frac{i}{r}w_{\theta}(z)\right),\ \ \  \overline{w_{\bar{z}}(z)}=\frac{e^{-i\theta}}{2}\left(\overline{w_r(z)}-\frac{i}{r}\overline{w_{\theta}(z)}\right)\ee
are analytic in $\ID$.

For each $\alpha\in[0,2\pi]$, the {\it directional derivative} of $w$ at $z$ is defined by
\beqq\label{eq1.2}
\partial_{\alpha}w(z)=\lim\limits_{r\rightarrow 0^+}\frac{w(z+re^{i\alpha})-w(z)}{re^{i\alpha}}=w_z(z)+e^{-2i\alpha}w_{\bar{z}}(z).
\eeqq
Then
\beqq\label{eq1.3}
\Lambda_w(z):=\max\limits_{0\leq\alpha\leq2\pi}\{|\partial_{\alpha}w(z)|\}=|w_z(z)|+|w_{\bar{z}}(z)|
\eeqq
and
\beqq\label{eq1.4}
\lambda_w(z):=\min\limits_{0\leq\alpha\leq2\pi}\{|\partial_{\alpha}w(z)|\}=\big||w_z(z)|-|w_{\bar{z}}(z)|\big|.
\eeqq
It is well known that $w$ is locally univalent and sense-preserving in $\ID$ if and only if its Jacobian satisfies
$$J_w(z)=|w_z(z)|^2-|w_{\bar{z}}(z)|^2>0,\ \ \  ~\mbox{for any}~\  z \in \ID.$$

\subsection*{Quasiregular mappings}
In order to state our motivations and results more precisely, we should introduce the definition of $n$-dimensional {\it quasiregular mappings}.
Following the definition in \cite[Page 127]{Matti} (see also \cite[Page 11 and Page 48]{Vaisala}), the definition of a quasiregular mapping in
a domain of $\mathbb{R}^n$ is given as follows:

Let $G\subset \mathbb{R}^n$ be a domain, and let $n\geq2$. A mapping $f:G\to \mathbb{R}^n$ is said to be {\it quasiregular} (briefly, qr.) if
\begin{enumerate}
  \item[(i)] $f$ is an {\it absolutely continuous} function in every line segment parallel to the coordinate axis and there exists the partial derivatives which are locally $L^n$ integrable functions on $\Omega$\  (we write $f\in ACL^n$).
  \item[(ii)] there exists a constant $K\geq1$ such that
  \be\label{May-4-1}L_f(x)^n\leq K J_f(x),\ee
  a.e. in $G$, where $L_f(x)$ is the {\it maximum stretching} for $f$ at the point $x$, i.e.,
  $$L_f(x)=\dlimsup\limits_{y\rightarrow x}\frac{|f(y)-f(x)|}{|y-x|},$$
and $J_f$ denotes the Jacobian determinant.
\end{enumerate}
If further, $f$ is a {\it homeomorphism} in $G$, then $f$ is said to be {\it quasiconformal}.

The smallest constant $K\geq1$ for which (\ref{May-4-1}) holds true is called the {\it outer dilatation} of $f$ and denote by $K_O(f)$.
If $f$ is quasiregular, then the smallest constant $K\geq1$, for which the inequality
$$J_f(x)\leq K l_f(x)^n,\ \ \ \mbox{where}\ \ \ l_f(x)=\min\{|f'(x)h|: |h|=1\},$$
holds a.e. in $G$, is called the {\it inner dilatation} of $f$ and denoted by $K_I(f)$. The {\it maximal dilatation} of $f$ is the number
$K(f)=\max\{K_I(f), K_{O}(f)\}$. If $K(f)\leq K$, then $f$ is said to be $K$-quasiregular ($K$-qr.). If $f$ is not quasiregular, we set $K_O(f)=K_I(f)=K(f)=\infty$.

It should be noted that the condition $f\in ACL^n$ guarantees the existence of the first derivatives of $f$ almost everywhere.
Moreover, the condition (i) is equivalent with the fact that $f$ is continuous and belongs to the Sobolev space $W^{1, n}_{loc}(G)$, i.e.,
the weak derivative is locally $L^n$ integrable in $G$, see for example \cite[Page 24 and Page 77]{Astala}.

Harmonic mappings and quasiconformal mappings are natural generalizations of conformal mappings. Harmonic mappings have nice algebraic properties
like power series and Poisson representation while quasiconformal mappings allows composition of mappings.
We refer the interested readers to \cite{Matti} for more discussions on
the conformal invariant of quasiregular mappings, and we refer to \cite{Kalaj-Jordan-08,DKal2011TAMS,Sakan-07,Pavlovic-02} for more discussions on harmonic quasiconformal mappings.

\subsection*{Motivations}
It was proved in \cite[Lemma 2.1]{Sakan-07} that if $w$ is a harmonic quasiconformal mapping of $\ID$ onto $\Omega\subset\mathbb{C}$, where
$\Omega$ is bounded by a rectifiable Jordan curve $\Gamma$, then $w_z\in H^1$ and $\overline{w_{\bar{z}}}\in H^1$.

Gehring showed in \cite[Theorem 1]{Gehring} that suppose $E$ is a domain in $\mathbb{R}^n$ and that $f:E\rightarrow \mathbb{R}^n$ is
a $K$-quasiconformal mapping. Then its maximum stretching $L_f$ is locally $L^p$-integrable in $E$ for $p\in[n, n+c)$, where $n\geq2$ and $c$ is a positive constant which depends only on $K$ and $n$.

Let
$$\mathrm{a_f}(x)=\mathrm{exp}\left[\frac{1}{n|B_x|}\int_{B_x}\log J_f(y)\mathrm{d}m(y)\right],$$
where $|B_x|$ is the $n$-measure of $B_x$. Notice that if $f$ is conformal, then the mean value property implies that $\mathrm{a_f}=L_f$.
It is easy to see that if $n=2$, then $L_f=\Lambda_f$.

Suppose $f$ is a quasiconformal mapping of $\mathbb{B}^n$ and fix $0<p<\infty$. Let
$F(\omega)=\lim_{r\to 1}f(r\omega)$ be the boundary function of $f$, where $\omega\in\mathbb{S}^{n-1}$, and set
$$\Gamma(\omega)=\{x\in\mathbb{B}^n: |x-\omega|\leq3(1-|x|)\}$$
be the cone with vertex $\omega$. Then, it follows from \cite[Theorem 5.1]{Astala-Hp} that the following conditions are equivalent:
(a) $F\in L^p(\mathbb{S}^{n-1})$; (b) $\int_{\mathbb{B}^n}\mathrm{a_f}(x)^p(1-|x|)^{p-1}\mathrm{d}m(x)<\infty$; (c) $\sup_{x\in\Gamma(\omega)}\mathrm{a_f}(x)(1-|x|)\,\in L^p(\mathbb{S}^{n-1})$.
Moreover, according to \cite[Theorem 9.3]{Astala-Hp}, we see that if $f\in L^{pn/(n-p)}(\mathbb{B}^n)$, $0<p<n$, then $L_f\in L^q(\mathbb{B}^n)$ for all $q<p$.
Finally, the authors in \cite{Astala-Hp} also presented three open problems related to quasiconformal mappings and the $H^p$ space.
We also refer to \cite{Koskela-19} for more discussions on weighted Hardy spaces and quasiconformal mappings.
It should be noted that in \cite{Astala-Hp, Koskela-19}, the condition $f$ is univalent, plays an important role in their proofs, see for example \cite[Lemma 2.1 and Lemma 2.3]{Astala-Hp} and \cite[Lemma 2.1 and Lemma 2.2]{Koskela-19}.

By comparing the above results, the following problem becomes interesting:

\bprob
Under what conditions on the boundary function $F$ ensure that the partial derivatives of its harmonic extension $w$,
i.e., $w_z$ and $\overline{w_{\bar{z}}}$, are in the space $L^p(\ID)$\, (or $H^p(\ID)$), where $p\geq1$?
\eprob

Suppose $w=P[F]$ is harmonic in $\ID$ with the boundary function $F$ is absolutely continuous.
Then, it follows from \cite[Chapter 6]{Rudin} that $F$ is a function of bounded variation. Thus, for
almost all $e^{it}\in \mathbb{T}$, the derivative  $\dot{F}(e^{it})$ exists, where
$$\dot{F}(e^{it}):=\frac{\mathrm{d}}{\mathrm{d}t}F(e^{it}).$$
Furthermore, we assume that $\dot{F}$ is of $L^p(0, 2\pi)$ space\, ($p\geq1$).

In this paper, under these assumptions on $F$, we prove that both $w_z$ and
$\overline{w_{\bar{z}}}$ are of $L^p(\ID)$ space for any $1\leq p<2$. Furthermore, if $w$ is a harmonic quasiregular mapping,
we show that both $w_z$ and $\overline{w_{\bar{z}}}$ are of $H^p$ space, for all $1\leq p\leq \infty$.
The Bergman norm estimates: $\|w_z\|_{L^p}$, $\|\overline{w_{\bar{z}}}\|_{L^p}$, and the Hardy norm estimates: $\|w_z\|_{p}$, $\|\overline{w_{\bar{z}}}\|_{p}$
are also obtained.
The main technique of this paper is
the Poisson integral, and in our proof, we do not require that $w$ is univalent.

Our main results are as follows:

\begin{thm}\label{thm1-2019-July-12}
Suppose $1\leq p<\infty$, $w=P[F]$ is a harmonic mapping of $\ID$ with the boundary function $F$ is absolutely continuous and satisfies $\dot{F}\in L^p(0, 2\pi)$.
Then for $z=re^{i\theta}\in\ID$,
$$\|w_r\|_{L^p}\leq (2C(p))^{1/p}\|\dot{F}\|_{L^p},$$
where $C(p)$ is a function of $p$ which is given by $($\ref{cp}$)$, and thus, $w_r(z)\in L^p(\ID).$
\end{thm}

\begin{rem}
(1) In Theorem $\ref{thm1-2019-July-12}$, the condition: ``$F$ is absolutely continuous" can not be weakened as: ``$F$ is of bounded variation".
This can be seen as follows:
If a function $F$ is of bounded variation, then $F$ has the following representation: $F=F_1+F_2$,
where $F_1$ is absolutely continuous and $F_2$ is completely singular, i.e., $\dot{F}_2=0$ a.e. on $\mathbb{T}$ (cf. \cite[Chapter 6]{Rudin}).
Now, suppose $F$ is completely singular. Then $\|\dot{F}\|_{L^p}=0$ a.e.
This implies that $w_r=w_{\theta}=0$, and thus, $w$ is a constant function.
However, there exists a function with its boundary function $F$ is
completely singular but its Poisson extension $P[F]$ is not a constant function (cf. \cite[Pages 58-62]{Du-04}).
Therefore, we should assume $F$ is absolutely continuous, which excludes the case of $F$ is completely singular.

(2) For the case $p=\infty$, the condition $\dot{F}\in L^{\infty}(0, 2\pi)$ can not ensure $w_r\in L^\infty(\ID)$.
This can be seen as follows: Suppose $F(e^{it})=|\sin t|$, where $t\in [0, 2\pi]$. Then $\dot{F}(t)=\cos t$, a.e. in $[0, 2\pi]$,
which shows that $\dot{F}\in L^{\infty}(0, 2\pi).$
However, elementary calculations show that
$$w=P[F](r)=\frac{1-r^2}{\pi r}\log\frac{1+r}{1-r}, \ \ \ 0<r<1.$$
Thus
$$w_r=\frac{2r-(1+r^2)\log\frac{1+r}{1-r}}{\pi r^2}\to \infty,$$
as $r\to 1$.

Moreover, this example also shows that $rw_r\notin h^p$, for any $1\leq p\leq\infty$.
\end{rem}

\begin{thm}\label{thm2-2019-July-12}
Suppose $1\leq p<2$, $w=P[F]$ is a harmonic mapping of $\ID$ with the boundary function $F$ is absolutely continuous and satisfies  $\dot{F}\in L^p(0, 2\pi)$.
Then
$$\|w_z\|_{L^p}\leq \left(C(p)+\frac{1}{2-p}\right)^{1/p}\|\dot{F}\|_{L^p}\ \ \mbox{and}\ \  \|\overline{w_{\bar{z}}}\|_{L^p}\leq \left(C(p)+\frac{1}{2-p}\right)^{1/p}\|\dot{F}\|_{L^p}$$
where $C(p)$ is given by $($\ref{cp}$)$, and this shows that $w_z, \overline{w_{\bar{z}}}\in L^p(\ID).$
\end{thm}

\begin{thm}\label{thm3-2019-July-14}
Suppose $1\leq p\leq\infty$, $w=P[F]$ is a harmonic quasiregular mapping of $\ID$ with the boundary function $F$ is absolutely continuous and satisfies $\dot{F}\in L^p(0, 2\pi)$.
Then
$$\|w_z\|_p\leq K\|\dot{F}\|_{L^p}\ \ \ \mbox{and}\ \ \ \|\overline{w_{\bar{z}}}\|_p\leq \frac{K-1}{2}\|\dot{F}\|_{L^p},$$
where $K\geq1$ is the outer dilatation of $w$. This shows that $ w_z\in H^p$ and $\overline{w_{\bar{z}}}\in H^p.$
\end{thm}

\begin{rem}
In Theorem \ref{thm3-2019-July-14}, the assumption that $w=P[F]$ is quasiregular can not be removed.
We use an example (Example \ref{exam4.1}, see also \cite[Page 62]{Sakan-97})
in Section \ref{sec-4} to show that there exists an absolutely continuous function $F$ satisfying $\dot{F}\in L^\infty(0, 2\pi)$ and $w=P[F]$ is harmonic in $\ID$ but not
quasiregular in $\ID$, and $w_z\notin L^{\infty}(\ID)$.
\end{rem}

\section{Preliminaries}\label{sec-2}
In this section, we should recall some known results and prove three lemmas.
We begin with the convex functions and Jensen's inequality.

\begin{defn}$($\cite[Definition 1]{Mitrinovic}$)$
(a) Let $I$ be an interval in $\mathbb{R}$. Then $f: I\rightarrow\mathbb{R}$ is said to be {\it convex}
if for all $x, y\in I$ and
$\lambda\in[0, 1]$,
\be\label{Junly-16-1}
f(\lambda x+(1-\lambda)y)\leq\lambda f(x)+(1-\lambda)f(y).
\ee
If (\ref{Junly-16-1}) is strict for all $x\neq y$ and $\lambda\in(0, 1)$, then $f$ is said to be {\it strictly convex}.

(b) If the inequality in (\ref{Junly-16-1}) is reversed, then $f$ is said to be {\it concave}. If it is strict for all $x\neq y$ and $\lambda\in(0, 1)$, then $f$ is said to be {\it strictly concave}.
\end{defn}

For $1\leq p<\infty$, the function $f(x)=x^p$ is convex in $(0, \infty)$. Thus,
for any $a, b>0$, the following inequality holds
\be\label{convex-ineq}
\left(\frac{a+b}{2}\right)^p\leq\frac{a^p+b^p}{2}.
\ee

\noindent{\bf Jensen's inequality} (See \cite{Fink} and \cite{Mitrinovic}).
Suppose $\mu$ is a regular Borel measure such that $\int_a^b\mathrm{d}\mu>0$, $f\in L^1(\mathrm{d}\mu)$, ie., $\int_a^bf(x)\mathrm{d}\mu$ exists, $\varphi$ is
a convex function. Then
$$\varphi\left(\frac{\int_a^bf(x)\mathrm{d}\mu}{\int_a^b\mathrm{d}\mu}\right)\leq\int_a^b\varphi(f(x))\mathrm{d}\mu\bigg/\int_a^b\mathrm{d}\mu.$$

Jensen's inequality has many applications. For example, assume that $f$ is a p.d.f. ({\it probability density function}) of a real-valued random
variable $X$, i.e., $f(x)\geq0$ and
$$\int_{-\infty}^{\infty}f(x)\mathrm{d}x=1,$$
$g$ is a continuous function and $\varphi$ is a convex function. Then
\be\label{Jensen-ineq}\varphi\left(\int_{-\infty}^{\infty}g(x)f(x)\mathrm{d}x\right)\leq\int_{-\infty}^{\infty}\varphi(g(x))f(x)\mathrm{d}x.\ee
This shows that
$$\varphi(E[g(X)])\leq E[\varphi\circ g(X)],$$
where $E(X)$ is the expectation of the random variable $X$.

\noindent{\bf Inverse hyperbolic tangent function.}
The function
$$\tanh x=\frac{e^x-e^{-x}}{e^x+e^{-x}}$$
is called the {\it hyperbolic tangent function}.
It is easy to see that $\tanh x$ is an odd, increasing function.
The Taylor series of $\tanh x$ is as follows:
$$\tanh x=\sum\limits_{n=1}^{\infty}\frac{2^{2n}(2^{2n}-1)B_{2n}x^{2n-1}}{(2n)!}=x-\frac{x^3}{3}+\frac{2x^5}{15}-\frac{17x^7}{315}+\cdots,\ \ \ \ |x|<\frac{\pi}{2},$$
where $B_{m}$ is the Bernoulli number which is defined by the following equation:
$$\frac{z}{e^{z}-1}=\sum\limits_{m=0}^{\infty}B_m\frac{z^m}{m!},\ \ \ z\in\mathbb{C}.$$
For some $m$, we can list the values of $B_m$ as follows: $B_0=1$, $B_1=-\frac{1}{2}$, $B_2=\frac{1}{6}$, $B_4=-\frac{1}{30}$, $B_6=\frac{1}{42}$, $\cdots$.
Moreover, $B_{2k+1}=0$, where $k\geq1$ is an integer.

The {\it inverse hyperbolic tangent function} is as follows:
$$\tanh^{-1} x=\frac{1}{2}\log\frac{1+x}{1-x}.$$
It is easy to see that
$\frac{\mathrm{d }}{\mathrm{\mathrm{d}}x}\tanh^{-1}x=\frac{1}{1-x^2}$ and $\tanh^{-1} x$ has the following Taylor series
$$\tanh^{-1}x=\sum\limits_{n=0}^{\infty}\frac{1}{2n+1}x^{2n+1}=x+\frac{x^3}{3}+\frac{x^5}{5}+\frac{x^7}{7}+\cdots,\ \ \ \ |x|<1.$$

\begin{lem}\label{lem1-2019-July-16}
For $0< r<1$, let
$$\varphi(r)=\log\frac{1}{1-r}-\frac{2\tanh^{-1}r}{r}.$$
Then $\varphi(r)$ is an increasing function of $r$.
\end{lem}
\bpf
Elementary calculations show that
$$\varphi'(r)=\frac{-r(2+r)+2(1+r)\tanh^{-1}r}{r^2(1+r)}.$$
The function $\psi(r):=-r(2+r)+2(1+r)\tanh^{-1}r$ is an increasing function of $r\in(0, 1)$, since
$$\psi'(r)=\frac{2[r^2+(1-r)\tanh^{-1}r]}{1-r}>0.$$
Therefore $\psi(r)> \psi(0)=0$, which shows that $\varphi'(r)>0$ for any $0<r<1$.

The proof of Lemma \ref{lem1-2019-July-16} is complete.
\epf

\begin{lem}\label{lem2-2019-July-12}
For $1\leq p<\infty$, $\theta\in[0, 2\pi]$ and $0\leq r<1$, let
$$I(r)=\frac{1}{\pi}\int_0^{2\pi}\frac{|\sin(t-\theta)|}{1+r^2-2r\cos(t-\theta)}\mathrm{d}t.$$
Then
\be\label{Ir}I(r)=\frac{4\tanh^{-1}r}{\pi r}\ee
and
\be\label{cp}C(p):=\int_0^1 I(r)^pr\mathrm{d}r\leq \frac{4^{p-1}}{\pi^p}\Big[2^{p}+(2-2^{-p})\Gamma(1+p)\Big].\ee
\end{lem}
\bpf
Elementary calculations show that
$$I(r)=\frac{2}{\pi}\int_0^{\pi}\frac{\sin x}{1+r^2-2r\cos x}\mathrm{d}x=\frac{2}{\pi r}\log\frac{1+r}{1-r},$$
and thus,
$$\int_0^1 I(r)^pr\mathrm{d}r=\left(\frac{2}{\pi}\right)^p\int_0^1\left(\frac{2 \tanh^{-1}r}{r}\right)^pr\mathrm{d}r.$$
It follows from Lemma \ref{lem1-2019-July-16} that
$$\varphi(r)=\log\frac{1}{1-r}-\frac{2\tanh^{-1}r}{r}$$
is an increasing function of $r\in[0, 1]$.
Therefore, $\varphi(r)\geq\varphi(0)=-2$, that is,
\be\label{July-13}\frac{2 \tanh^{-1}r}{r}\leq 2+\log\frac{1}{1-r}.\ee
For $p\geq1$, using (\ref{convex-ineq}) we have the following inequality
\be\label{July-13-1}\left(\frac{2+\log\frac{1}{1-r}}{2}\right)^p\leq\frac{2^p+\left(\log\frac{1}{1-r}\right)^p}{2}.\ee
Then the inequalities (\ref{July-13}) and (\ref{July-13-1}) lead to
$$\int_0^1 I(r)^pr\mathrm{d}r\leq\left(\frac{2}{\pi}\right)^p\int_0^12^{p-1}\left[2^p+\left(\log\frac{1}{1-r}\right)^p\right]r\mathrm{d}r.$$
Recall that for $p\geq1$ and $\alpha>-1$, the following equality holds
$$\int_0^1 t^{\alpha}\left(\log\frac{1}{t}\right)^{p-1}\mathrm{d}t=\frac{\Gamma(p)}{(1+\alpha)^p}.$$
Then
$$\int_0^1 \left(\log\frac{1}{1-r}\right)^pr\mathrm{d}r=\frac{2-2^{-p}}{2}\Gamma(1+p).$$
Based on the above facts, we have
$$\int_0^1 I(r)^pr\mathrm{d}r\leq\frac{4^{p-1}}{\pi^p}\Big[2^{p}+(2-2^{-p})\Gamma(1+p)\Big].$$

This completes the proof of Lemma \ref{lem2-2019-July-12}.
\epf

For some positive integers $p$, we list some values of the function $C(p)$ as follows:
\begin{table}[h]
\centerline{\large\begin{tabular}{|c|c|c|c|c|c|c|}
\hline
\ $p$\ &\ 1\ &\ 2\ &\ 3\ &\ 4\ &\ 5\ \\
\hline
\ $C(p)$\ &\ $\frac{\pi}{2}$ \ &\ $\frac{8}{3}$\ &\  $\frac{16}{\pi}$\ &\ $\frac{128(30+\pi^2)}{45\pi^2}$\ &\ $\frac{256(15+2\pi^2)}{9\pi^2}$\  \\
\hline
\end{tabular}}
\end{table}

\begin{lem}\label{thm0-2019-July-12}
Suppose $1\leq p\leq\infty$, $w=P[F]$ is a harmonic mapping of $\ID$ with the boundary function $F$ is absolutely continuous and satisfies $\dot{F}\in L^p(0, 2\pi)$.
Then for $z=re^{i\theta}\in\ID$,
$$\|w_{\theta}\|_p\leq \|\dot{F}\|_{L^p},$$
and thus, $w_{\theta}(z)\in h^p.$
\end{lem}
\bpf
For $z=re^{i\theta}\in\ID$, integral by part leads to
\begin{align*}
w_{\theta}(re^{i\theta})&=\int_0^{2\pi}\frac{\partial}{\partial \theta}\left\{P_r(t-\theta)\right\}F(e^{it})\mathrm{d}t\\
&=-\int_0^{2\pi}F(e^{it})\frac{\partial}{\partial t}\left\{P_r(t-\theta)\right\}\mathrm{d}t\\
&=\int_0^{2\pi}P_r(t-\theta)\mathrm{d}F(e^{it}).
\end{align*}
By using $\int_0^{2\pi}P_r(t-\theta)\mathrm{d}t=1$ and Jensen's inequality (note that for $1\leq p<\infty$, $\varphi(x)=x^p$ is convex), we have
$$\left|w_{\theta}(re^{i\theta})\right|^p\leq\left(\int_0^{2\pi}P_r(t-\theta)|\dot{F}(e^{it})|\mathrm{d}t\right)^p\leq\int_0^{2\pi}P_r(t-\theta)|\dot{F}(e^{it})|^p\mathrm{d}t,$$
where $1\leq p<\infty$.
The assumption of $\dot{F}\in L^p(0, 2\pi)$ ensures that
$$P_r(t-\theta)|\dot{F}(e^{it})|^p\in L^1((0, 2\pi)\times (0, 2\pi)).$$
Using Fubini's Theorem we obtain that
\begin{align}\label{LpF}
\int_{0}^{2\pi}\left|w_{\theta}(re^{i\theta})\right|^p\mathrm{d}\theta &\leq\int_0^{2\pi}\mathrm{d}\theta\int_0^{2\pi}P_r(t-\theta)|\dot{F}(e^{it})|^p\mathrm{d}t\\ \nonumber
&=\int_0^{2\pi}|\dot{F}(e^{it})|^p\mathrm{d}t\int_0^{2\pi}P_r(t-\theta)\mathrm{d}\theta\\ \nonumber
&=2\pi\|\dot{F}\|_{L^p}^p,
\end{align}
which shows that
$$\|w_{\theta}\|_p=\lim_{r\to 1}\left(\frac{1}{2\pi}\int_0^{2\pi}\left|w_{\theta}(re^{i\theta})\right|^p\mathrm{d}\theta\right)^{1/p}\leq \|\dot{F}\|_{L^p},$$
and thus, $w_\theta(z)\in h^p$.

For the case $p=\infty$, since $F$ is absolutely continuous, we see from \cite[Page 100]{Pavlovic-02} that
$\lim_{r\rightarrow1}w_{\theta}(re^{i\theta})=\dot{F}(e^{i\theta})$ a.e. on $[0, 2\pi]$.
Note that $w_{\theta}$ is a harmonic mapping of $\ID$, then the {\it maximum principle} property shows that
\be\label{may-6}\|w_{\theta}\|_{\infty}\leq\|\dot{F}\|_{\infty},\ee
which proves $w_{\theta}\in h^{\infty}$.

The proof of Lemma \ref{thm0-2019-July-12} is complete.
\epf

Let us end this section by recalling the following results which show that it is natural to assume the boundary function $F$ is
absolutely continuous when we consider harmonic quasiregular mappings of $\mathbb{D}$ onto a bounded domain $\Omega\subset\mathbb{C}$.

Recall that the Cauchy singular integral $C_{\mathbb{T}}[\varphi]$ of a function $\varphi:\mathbb{T}\rightarrow\mathbb{C}$, which is Lebesgue integrable
on $\mathbb{T}$, is defined as follows: for every $\zeta\in\mathbb{T}$, let
\be\label{CT}C_{\mathbb{T}}[\varphi](\zeta):=p.v.\frac{1}{2\pi i}\int_{\mathbb{T}}\frac{\varphi(u)}{u-\zeta}\mathrm{d}u:=\lim\limits_{\epsilon\rightarrow0^+}\frac{1}{2\pi i}\int_{\mathbb{T}\setminus \mathbb{T}(\zeta, \epsilon)}\frac{\varphi(u)}{u-\zeta}\mathrm{d}u\ee
whenever the limit exists, and $C_{\mathbb{T}}[\varphi](\zeta):=0$ otherwise, where $ \mathbb{T}(e^{ix}, \epsilon):=\{e^{it}\in\mathbb{T}: |t-x|<\epsilon\}$.

Given a continuous function $\varphi:\mathbb{T}\rightarrow \mathbb{C}$ and $\zeta\in\mathbb{T}$, set
\be\label{Vf}V[\varphi](\zeta):=\lim\limits_{\epsilon\rightarrow0^+}\frac{1}{2\pi}\int_{\mathbb{T}\setminus \mathbb{T}(\zeta, \epsilon)}\frac{|\varphi(u)-\varphi(\zeta)|^2}{|u-\zeta|^2}|\mathrm{d}u|,\ee
and
\be\label{Vs}V^*[\varphi](\zeta):=-\lim\limits_{\epsilon\rightarrow0^+}\frac{1}{\pi}\int_{\mathbb{T}\setminus \mathbb{T}(\zeta, \epsilon)}\frac{\textbf{Im}(\varphi(u)\overline{\varphi(\zeta)})}{|u-\zeta|^2}|\mathrm{d}u|,\ee
provided the limits exist, as well as $V[\varphi](\zeta):=\infty$ and $V^*[\varphi](\zeta):=0$ otherwise.

\vspace*{3mm}
\noindent{\bf Theorem A. }$($\cite[Theorem 1.2]{Sakan-07}$)$
{\it If $F$ is a homeomorphism of $\mathbb{T}$ and absolutely continuous on $\mathbb{T}$, then for a.e.
$\zeta\in\mathbb{T}$ the limit in $($\ref{CT}$)$ with $\varphi$ replaced by $F'$ and the limits in $($\ref{Vf}$)$ and $($\ref{Vs}$)$
exist, and}
$$2C_{\mathbb{T}}[F'](\zeta)=\bar{\zeta}F(\zeta)(V[F](\zeta)+iV^*[F](\zeta)).$$

\vspace*{3mm}
\noindent{\bf Corollary B. }$($\cite[Corollary 2.2]{Sakan-07}$)$
{\it Given $K\geq1$ and a domain $\Omega$ in $\mathbb{C}$, let $w=P[F]$ be a harmonic quasiconformal mapping of $\ID$ onto $\Omega$.
If $\Omega$ is bounded by a rectifiable Jordan curve $\Gamma$, then $F$ is absolutely continuous.}

\section{Proofs of the main results}\label{sec-3}

\subsection*{Proof of Theorem \ref{thm1-2019-July-12}}
Since $F$ is absolutely continuous, integral by part shows that
\begin{align*}
w_r(re^{i\theta})&=\int_0^{2\pi}\frac{\partial}{\partial r}\{P_r(t-\theta)\}F(e^{it})\mathrm{d}t\\ \nonumber
&=\frac{2}{1-r^2}\int_0^{2\pi}\frac{\partial}{\partial t}\{P_r(t-\theta)\sin(t-\theta)\}F(e^{it})\mathrm{d}t\\ \nonumber
&=\frac{2}{r^2-1}\int_0^{2\pi}P_r(t-\theta)\sin(t-\theta)\dot{F}(e^{it})\mathrm{d}t.
\end{align*}
Thus
\begin{align*}
\left|w_r(re^{i\theta})\right|&\leq\frac{2}{1-r^2}\int_0^{2\pi}P_r(t-\theta)|\sin(t-\theta)||\dot{F}(e^{it})|\mathrm{d}t\\ \nonumber
&=\frac{1}{\pi}\int_0^{2\pi}\frac{|\sin(t-\theta)|}{1+r^2-2r\cos(t-\theta)}|\dot{F}(e^{it})|\mathrm{d}t.
\end{align*}
Let
$$I(r)=\frac{1}{\pi}\int_0^{2\pi}\frac{|\sin(t-\theta)|}{1+r^2-2r\cos(t-\theta)}\mathrm{d}t.$$
It follows from (\ref{Ir}) that
$$I(r)=\frac{4\tanh^{-1}r}{\pi r}.$$
For $1\leq p<\infty$, according to Jensen's inequality (note that $\varphi(x)=x^p$ is convex), we have
\begin{align*}
\left|w_r(re^{i\theta})\right|^p&\leq  I(r)^p\left(\frac{1}{\pi}\int_0^{2\pi}\frac{|\sin(t-\theta)|}{1+r^2-2r\cos(t-\theta)}\frac{1}{I(r)}|\dot{F}(e^{it})|\mathrm{d}t\right)^p\\ \nonumber
&\leq \frac{I(r)^{p-1}}{\pi}\int_0^{2\pi}\frac{|\sin(t-\theta)|}{1+r^2-2r\cos(t-\theta)}|\dot{F}(e^{it})|^p\mathrm{d}t.
\end{align*}
The assumption of $\dot{F}\in L^p(0, 2\pi)$\,  ($1\leq p<\infty$) ensures that
$$\frac{|\sin(t-\theta)|}{1+r^2-2r\cos(t-\theta)}|\dot{F}(e^{it})|^p\in L^1((0, 2\pi)\times (0, 2\pi)).$$
By using Fubini's Theorem, we obtain that
\begin{align*}
\int_{0}^{2\pi}\left|w_r(re^{i\theta})\right|^p\mathrm{d}\theta &\leq\frac{I(r)^{p-1}}{\pi}\int_0^{2\pi}\mathrm{d}\theta\int_0^{2\pi}\frac{|\sin(t-\theta)|}{1+r^2-2r\cos(t-\theta)}|\dot{F}(e^{it})|^p\mathrm{d}t\\ \nonumber
&=\frac{I(r)^{p-1}}{\pi}\int_0^{2\pi}|\dot{F}(e^{it})|^p\mathrm{d}t\int_0^{2\pi}\frac{|\sin(t-\theta)|}{1+r^2-2r\cos(t-\theta)}\mathrm{d}\theta\\ \nonumber
&\leq2\pi\|\dot{F}\|_{L^p}^p I(r)^p,
\end{align*}
and thus,
\begin{align*}
  \int_{\ID}\left|w_r(re^{i\theta})\right|^p\mathrm{d}A(z)& =\frac{1}{\pi}\int_0^1r\mathrm{d}r\int_{0}^{2\pi}\left|w_r(re^{i\theta})\right|^p\mathrm{d}\theta  \\
   & \leq 2 \|\dot{F}\|_{L^p}^p \int_0^1 I(r)^pr\mathrm{d}r= 2 \|\dot{F}\|_{L^p}^p C(p),
\end{align*}
where $C(p)$ is given by (\ref{cp}).
Then
\be\label{July-13-2}\|w_r\|_{L^p}^p=\int_{\ID}\left|w_r(re^{i\theta})\right|^p\mathrm{d}A(z)\leq 2C(p)\|\dot{F}\|_{L^p}^p,\ee
which shows that
$$\|w_r\|_{L^p}\leq (2C(p))^{1/p}\|\dot{F}\|_{L^p},$$
and thus, $w_r(re^{i\theta})\in L^p(\ID)$.

The proof of Theorem \ref{thm1-2019-July-12} is complete.
\qed

\subsection*{Proof of Theorem \ref{thm2-2019-July-12}}
For $z=re^{i\theta}\in\ID$, it follows from (\ref{Oct-20-1}) that
$$ |w_z(z)|\leq\frac{1}{2}\left(|w_r(z)|+\left|\frac{w_{\theta}(z)}{r}\right|\right).$$
For $1\leq p<\infty$, applying (\ref{convex-ineq}) we have
$$|w_z(z)|^p\leq\frac{1}{2^p}\left(|w_r(z)|+\left|\frac{w_{\theta}(z)}{r}\right|\right)^p\leq\frac{1}{2}\left(|w_r(z)|^p+\left|\frac{w_{\theta}(z)}{r}\right|^p\right).$$
We first estimate
$$\int_{\ID}\left|\frac{w_{\theta}(z)}{r}\right|^p\mathrm{d}A(z)$$
as follows:
According to (\ref{LpF}), we see that
$$\int_0^{2\pi}|w_\theta(re^{i\theta})|^p\mathrm{d}\theta\leq2\pi\|\dot{F}\|_{L^p}^p.$$
This implies that
$$\int_{\ID}\left|\frac{w_{\theta}(re^{i\theta})}{r}\right|^p\mathrm{d}A(z)\leq2\|\dot{F}\|_{L^p}^p\int_0^1r^{1-p}\mathrm{d}r=\frac{2\|\dot{F}\|_{L^p}^p}{2-p},$$
where $1\leq p<2.$

On the other hand, we already showed in (\ref{July-13-2}) that
$$\int_{\ID}\left|w_r(re^{i\theta})\right|^p\mathrm{d}A(z)\leq 2 C(p) \|\dot{F}\|_{L^p}^p,$$
where $C(p)$ is given by (\ref{cp}) and
$$C(p)\leq \frac{4^{p-1}}{\pi^p}\Big[2^{p}+(2-2^{-p})\Gamma(1+p)\Big].$$
Based on these facts, we have
$$\int_{\ID}|w_z(z)|^p\mathrm{d}A(z)\leq\left(C(p)+\frac{1}{2-p}\right) \|\dot{F}\|_{L^p}^p,$$
which shows that
$$\|w_z\|_{L^p}\leq \left(C(p)+\frac{1}{2-p}\right)^{1/p}\|\dot{F}\|_{L^p},$$
and thus, $w_z\in L^p(\ID)$, for $1\leq p<2$.

Similarly, we can prove
$$\|\overline{w_{\bar{z}}}\|_{L^p}\leq \left(C(p)+\frac{1}{2-p}\right)^{1/p}\|\dot{F}\|_{L^p},$$
and thus, $\overline{w_{\bar{z}}}\in L^p(\ID)$, for $1\leq p<2$.

The proof of Theorem \ref{thm2-2019-July-12} is complete.
\qed

\subsection*{Proof of Theorem \ref{thm3-2019-July-14}}
For $z=re^{i\theta}\in\ID$, it follows from (\ref{5.3}) and (\ref{LpF}) that
$$\int_{0}^{2\pi}\left(|w_z(re^{i\theta})|-|w_{\bar{z}}(re^{i\theta})|\right)^p\mathrm{d}\theta\leq\int_{0}^{2\pi}\left|w_{\theta}(re^{i\theta})\right|^p\mathrm{d}\theta\leq 2\pi\|\dot{F}\|_{L^p}^p,$$
where $1\leq p<\infty$.
Since $w$ is a quasiregular mapping, we see that there exists a constant $K\geq1$ (the outer dilatation of $w$), such that
$$|w_z(re^{i\theta})|+|w_{\bar{z}}(re^{i\theta})|\leq K(|w_z(re^{i\theta})|-|w_{\bar{z}}(re^{i\theta})|).$$
Therefore,
$$\int_{0}^{2\pi}\left(|w_z(re^{i\theta})|+|w_{\bar{z}}(re^{i\theta})|\right)^p\mathrm{d}\theta\leq 2\pi K^p\|\dot{F}\|_{L^p}^p,$$
and thus,
$$\frac{1}{2\pi}\int_{0}^{2\pi}|w_z(re^{i\theta})|^p\mathrm{d}\theta\leq K^p\|\dot{F}\|_{L^p}^p,\ \ \ \frac{1}{2\pi}\int_{0}^{2\pi}|w_{\bar{z}}(re^{i\theta})|^p\mathrm{d}\theta\leq  \left(\frac{K-1}{2}\right)^p\|\dot{F}\|_{L^p}^p.$$
This shows that
$$M_p(r, w_z)\leq K\|\dot{F}\|_{L^p}\ \ \ \mbox{and}\ \ \ M_p(r, \overline{w_{\bar{z}}})\leq\frac{K-1}{2}\|\dot{F}\|_{L^p}.$$
Therefore, letting $r$ tends to 1, we have
$$\|w_z\|_p\leq K\|\dot{F}\|_{L^p}\ \ \ \mbox{and}\ \ \ \|\overline{w_{\bar{z}}}\|_p\leq \frac{K-1}{2}\|\dot{F}\|_{L^p},$$
which guarantee that $w_z\in H^p$ and $\overline{w_{\bar{z}}}\in H^p$, where $1\leq p<\infty$.

For the case $p=\infty$, by using (\ref{5.3}) and (\ref{may-6}), we see that
$$|w_z(re^{i\theta})|-|w_{\bar{z}}(re^{i\theta})|\leq |w_{\theta}(re^{i\theta})|\leq\|w_{\theta}\|_{\infty}\leq\|\dot{F}\|_{\infty}.$$
The quasiregularity of $w$ ensures that, there exists a constant $K\geq1$, such that
$$|w_z(re^{i\theta})|+|w_{\bar{z}}(re^{i\theta})|\leq K\left(|w_z(re^{i\theta})|-|w_{\bar{z}}(re^{i\theta})|\right)\leq K\|\dot{F}\|_{\infty}.$$
Then $\|w_z\|_{\infty}\leq K\|\dot{F}\|_{\infty}$, and thus, $w_z\in H^{\infty}$.

Similarly, we can prove $\|\overline{w_{\bar{z}}}\|_{\infty}\leq \frac{K-1}{2}\|\dot{F}\|_{L^{\infty}}$, and thus,
$\overline{w_{\bar{z}}}\in H^{\infty}$.

The proof of Theorem \ref{thm3-2019-July-14} is complete.
\qed

\section{An example}\label{sec-4}
In the following, we are going to construct an example (cf. \cite[Page 62]{Sakan-97}), which shows that the condition $w$ is quasiregular
in Theorem \ref{thm3-2019-July-14} cannot be removed. Before we start our discussion, we need to do some preparations.

Following the notation in \cite{Pavlovic-book,Pavlovic-02}, suppose $\varphi$ is a continuous increasing function on $\mathbb{R}$, such that $\varphi(2\pi+x)-\varphi(x)\equiv2\pi$, and let $F$ be the boundary function on $\mathbb{T}$, satisfying
\be\label{May-4-boundary}F(e^{it})=\Phi(t)=e^{i\varphi(t)},\ee
where $\Phi$ is a $2\pi$-periodic, absolutely continuous function on $[0, 2\pi]$. According to \cite[Page 100]{Pavlovic-book}, we see that the {\it Hilbert transformation} of $\Phi'$, which is defined as follows (see for example \cite[(2.1)]{Pavlovic-02} or \cite[Page 242]{Kalaj-Jordan-08}):
$$H[\Phi'](\theta)=-\frac{1}{\pi}\lim\limits_{\epsilon\rightarrow 0^+}\int\limits_{\epsilon}^{\pi}\frac{\Phi'(\theta+t)-\Phi'(\theta-t)}{2\tan\frac{t}{2}}\,\mathrm{d}t,$$
exists almost everywhere. Moreover, we have $\lim_{r\to 1}rw_r(re^{i\theta})=H[\Phi'](\theta)$ a.e. on $[0, 2\pi]$, where $w(z)=P[F](z)$ and $z=re^{i\theta}\in\ID$.

Pavlovi\'c proved in \cite[Theorem 6.6.1]{Pavlovic-book} and \cite[Theorem 1.1]{Pavlovic-02} that the harmonic mapping $w=P[F]$, where $F$ is given by (\ref{May-4-boundary}), is quasiconformal if and only if $\Phi$ is absolutely continuous and satisfies the following conditions:
(i) $\essinf_{\theta\in[0, 2\pi]}|\Phi'(\theta)|>0$, (ii) $\esssup_{ \theta\in[0, 2\pi]} |\Phi'(\theta)|<\infty$, (iii) $\esssup_{\theta\in[0, 2\pi]}|H[\Phi'](\theta)|<\infty$.
Moreover, $w$ is quasiconformal if and only if $w$ is bi-Lipschitz.

Based on these results, we now use the following example to show that there exists a boundary function $F$, such that
$\dot{F}\in L^{\infty}(0, 2\pi)$, but $w=P[F]$ is not quasiregular (therefore, not quasiconformal), and
$w_z\notin L^{\infty}(\ID)$.

\begin{exam}\label{exam4.1}(\cite[Page 62]{Sakan-97})
Let
$$\varphi_0(x)=\left\{
\begin{array}
{r@{\ }l}
1+\left(1+\frac{1}{\pi}\right)x, \  & -\pi\leq x<0,\\
\\
1+\left(1-\frac{1}{\pi}\right)x,  \  & 0\leq x\leq\pi.
\end{array}\right.$$
For all $x\in [-\pi, \pi]$ and integer $k$, set $\varphi(x+2k\pi)=\varphi_0(x)+2k\pi$ and
$$F(e^{ix})=e^{i\varphi(x)}.$$
Then the function $\varphi:\mathbb{R}\rightarrow\mathbb{R}$ satisfies the following equation: $\varphi(x+2k\pi)=\varphi(x)+2k\pi$, where
$x\in\mathbb{R}$ and $k$ is an integer.
The following statements hold:
\ben
\item[(\rm{A1})]
$\dot{F}\in L^{\infty}(0, 2\pi)$;
\item[(\rm{A2})]
$w=P[F]$ is harmonic in $\ID$ but not quasiregular in $\ID$;
\item[(\rm{A3})]
$w_z\notin L^{\infty}(\ID)$.
\een

\end{exam}

\bpf
(A1)
Let $\zeta=e^{ix}$, where $x\in[-\pi, \pi]$. Then
$$F(\zeta)=e^{i\varphi(-i \log \zeta)}.$$  Here and hereafter, denote by $\log$ the principal value of the natural logarithm.
For any $x\neq0$, $\pm\pi$,
$$\left|\frac{\mathrm{d}}{\mathrm{d}x}F(e^{ix})\right|=|\varphi_0'(x)| \leq1+\frac{1}{\pi},$$
which shows $\dot{F}\in L^{\infty}(0, 2\pi)$.

(A2) Obviously, we have $w$ is a harmonic self-mapping of $\ID$. According to the definition of $\varphi$, we see that $\varphi$ is
an increasing, continuous function with its derivative exists a.e. on $\mathbb{R}$, and thus, $\varphi$ is absolutely continuous.

Now, we prove $w$ is not a quasiregular mapping by showing that the Hilbert transformation of $\dot{F}$ is not essentially bounded.

Let
$$\Phi(x)=e^{i\varphi(x)}.$$
Then $\Phi'(x)$ exists and continuous a.e. on $[-\pi, \pi]$.
Elementary calculations show that
\begin{eqnarray*}
|H[\Phi'](0)|&=&\lim\limits_{\epsilon\rightarrow0^+}\frac{1}{\pi}\left|\int\limits_{\epsilon}^{\pi}\frac{\Phi'(t)-\Phi'(-t)}{2\tan\frac{t}{2}}\,\mathrm{d}t\right|\\
&=&\lim\limits_{\epsilon\rightarrow0^+}\frac{1}{\pi}\left|\frac{1}{\pi}\int\limits_{\epsilon}^{\pi}\frac{e^{it}(e^{\frac{it}{\pi}}+e^{\frac{-it}{\pi}})}{2\tan\frac{t}{2}}\,\mathrm{d}t+\int\limits_{\epsilon}^{\pi}\frac{e^{it}(e^{\frac{it}{\pi}}-e^{\frac{-it}{\pi}})}{2\tan\frac{t}{2}}\,\mathrm{d}t\right| \\
&\geq&\frac{1}{\pi^2}\int_{\epsilon}^{\pi}\frac{\cos\frac{t}{\pi}}{\tan\frac{t}{2}}\,\mathrm{d}t-\frac{2}{\pi^2}\int_0^{\pi}\frac{\sin^2\frac{t}{2}\cos\frac{t}{\pi}}{\tan\frac{t}{2}}\mathrm{d}t\\
&&-\frac{1}{\pi^2}\int_0^{\pi}\frac{\sin t\cos\frac{t}{\pi}}{\tan\frac{t}{2}}\mathrm{d}t-\frac{1}{\pi}\int_{0}^{\pi}\frac{\sin\frac{t}{\pi}}{\tan\frac{t}{2}}\,\mathrm{d}t.
\end{eqnarray*}
It is easy to see that
$$\frac{2}{\pi^2}\int_0^{\pi}\frac{\sin^2\frac{t}{2}\cos\frac{t}{\pi}}{\tan\frac{t}{2}}\mathrm{d}t=\frac{1+\cos1}{\pi^2-1},$$
$$\frac{1}{\pi^2}\int_0^{\pi}\frac{\sin t\cos\frac{t}{\pi}}{\tan\frac{t}{2}}\mathrm{d}t=\frac{\sin1}{\pi^2-1},$$
and
$$\frac{1}{\pi}\int_{0}^{\pi}\frac{\sin\frac{t}{\pi}}{\tan\frac{t}{2}}\mathrm{d}t\leq\frac{2}{\pi}.$$
Then there is a constant $M=\frac{1+\cos1}{\pi^2-1}+\frac{\sin1}{\pi^2-1}+\frac{2}{\pi}>0$ such that
$$|H[\Phi'](0)|\geq \lim\limits_{\epsilon\rightarrow0^+}\frac{1}{\pi^2}\int_{\epsilon}^{\pi}\frac{\cos\frac{t}{\pi}}{\tan\frac{t}{2}}\,\mathrm{d}t-M.$$
The divergence of the integral $\int_{0}^{\pi}\frac{\cos\frac{t}{\pi}}{\tan\frac{t}{2}}\mathrm{d}t$ shows that
\be\label{Hilbert}|H[\Phi'](0)|=\infty.\ee
Since $H[\Phi'](x)$ continuous a.e. on $[-\pi, \pi]$,  we see that (cf. \cite[Page 62]{Sakan-97})
$$\esssup \{H[\Phi'](x): x\in[-\pi, \pi]\}=\infty.$$

Moreover, by straightforward computation we find that (cf. \cite[Page 100]{Pavlovic-book}) $|w_r(e^{i\theta})|^2=A(\theta)^2+B(\theta)^2$,
where
$$A(\theta)=\frac{1}{2\pi}\int_{-\pi}^{\pi}\left(\frac{\sin(\varphi(\theta+t)/2-\varphi(\theta)/2)}{\sin t/2}\right)^2\mathrm{d}t$$
and
$$B(\theta)=\frac{-1}{\pi}\int_{0^+}^{\pi}\frac{\sin(\varphi(\theta+t)-\varphi(\theta))+\sin(\varphi(\theta-t)-\varphi(\theta))}{4\sin^2(t/2)}\mathrm{d}t.$$
This implies that
$$\esssup_{\theta\in[-\pi, \pi]} (A(\theta)^2+B(\theta)^2)=\infty,$$
since $\lim_{r\to 1}|rw_r(re^{i\theta})|=|H[\Phi'](\theta)|$.

On the other hand, we already knew $|\varphi'(\theta)|<1+\frac{1}{\pi}$, and it follows from \cite[(6.26)]{Pavlovic-book} that
$$|w_z(e^{i\theta})|^2=\frac{1}{4}\bigg((A(\theta)+\varphi'(\theta))^2+B(\theta)^2\bigg)$$
and
$$|w_{\bar{z}}(e^{i\theta})|^2=\frac{1}{4}\bigg((A(\theta)-\varphi'(\theta))^2+B(\theta)^2\bigg).$$
Based on the above discussions, we have
$$\esssup_{z\in\ID}\left|\frac{w_{\bar{z}}(z)}{w_z(z)}\right|=1,$$
which shows that $w$ is not quasiregular.

(A3) As we have said before, $w$ is a quasiconformal self-mapping of $\ID$ if and only if $w$ is bi-Lipschitz.
In (A2), we already showed that $H[\Phi']$ is unbounded and $w$ is not quasiregular (and thus, not quasiconformal) in $\ID$.
Therefore, $w$ is not Lipschitz continuous in $\ID$, which implies that $w_z\notin L^{\infty}(\ID)$.

The proof of Example \ref{exam4.1} is complete.
\epf

\vspace*{5mm}
\noindent {\bf Acknowledgments}.
The author of this paper would like to thank the anonymous referee for his/her helpful comments that have significant impact on this paper, and would like to thank Professor Ken-ichi Sakan for his help on the discussions of Theorem \ref{thm1-2019-July-12}.

\vspace*{5mm}
\noindent{\bf Funding.}
The research of the author was supported by NSFs of China (No. 11501220), NSFs of Fujian Province (No. 2016J01020),
and the Promotion Program for Young and Middle-aged Teachers in Science and Technology Research of Huaqiao University (ZQN-PY402).


\begin{thebibliography}{1}
\bibitem{Astala}
K. Astala, T. Iwaniec, and G. Martin,  \emph{Elliptic partial differential equations and quasiconformal mappings in the plane},
Princeton Mathematical Series, Vol. 48, Princeton University Press, Princeton, NJ, 2009, p. xviii+677.

\bibitem{Astala-Hp}
K. Astala and P. Koskela, \emph{$H^p$-theory for quasiconformal mappings},
Pure Appl. Math. Q. \textbf{7} (2011), 19--50.

\bibitem{Koskela-19}
S. Benedict, P. Koskela, and X. Li, \emph{Weighted Hardy spaces of quasiconformal mappings},
arXiv:1904.00519, April, 2019.


\bibitem{Du-hp}P. Duren,
    \emph{Theory of $H^p$ spaces},
    Academic Press, New York, 1970.

\bibitem{Du-04} P. Duren,
\emph{Harmonic mappings in the plane},
{Cambridge Univ. Press}, New York, 2004.


\bibitem{Fink}
A. M. Fink and M. Jodeit, \emph{Jensen inequalities for functions with higher monotonicities},
Aequationes Math. \textbf{40} (1990), 26--43.

\bibitem{Gehring}
F. Gehring, \emph{The $L^p$-integrability of the partial derivatives of a quasiconformal mapping},
Acta Math. \textbf{130} (1973), 265--277.

\bibitem{Zhukh}
H. Hedenmalm, B. Korenblum, and K. Zhu,
\emph{Theory of Bergman spaces},
{Springer}, New York, 2000.

\bibitem{Kalaj-Jordan-08}D. Kalaj,
    \emph{ Quasiconformal and harmonic mappings between Jordan domains},
    Math. Zeit. \textbf{260} (2008), 237--252.

\bibitem{DKal2011TAMS}
D. Kalaj and M. Pavlovi\'c, \emph{{On quasiconformal self-mappings of the unit disk satisfying Poisson's equation}},
Trans. Amer. Math. Soc. \textbf{363} (2011), 4043--4061.

\bibitem{Mitrinovic}D. Mitrinovi\'c, J. Pe\v{c}ari\'c, and A. Fink,
    \emph{Classical and new inequalities in analysis},
    Kluwer Academic Publishers, London, 1993.

\bibitem{Sakan-97}D. Partyka and K. Sakan,
    \emph{ A note on non-quasiconformal harmonic extensions},
    Bull. Soc. Sci. Lettres Lodz, \textbf{47}(1997), 51--63.

\bibitem{Sakan-07}
D. Partyka and K. Sakan, \emph{On bi-Lipschitz type inequalities for quasiconformal harmonic mappings},
Ann. Acad. Sci. Fenn. Ser. A1-Math.  \textbf{32} (2007), 579--594.

\bibitem{Pavlovic-book} M. Pavlovi\'c,
 \emph{Introduction to function spaces on the disk},
Beograd Publisher, Beograd, 2004.


\bibitem{Pavlovic-02} M. Pavlovi\'c,
 \emph{Boundary correspondence under harmonic quasiconformal homeomorphisma of the unit disk},
Ann. Acad. Sci. Fenn. Ser. A1-Math.  \textbf{27} (2002), 365--372.

\bibitem{Rudin} W. Rudin,
\emph{Real and complex analysis},
{McGraw-Hill Education}, New York, 1986.

\bibitem{Vaisala} J. V\"ais\"al\"a,
\emph{Lectures on n-dimensional quasiconformal mappings},
Springer, Berlin Heidelberg, 1971.


\bibitem{Matti} M. Vuorinen,
\emph{Comformal geometry and quasiregular mappings},
Springer, Berlin Heidelberg, 1988.

\bibitem{Yamashita}
S. Yamashita, \emph{Hardy norm, Bergman norm, and univalency},
Ann. Polonici Mathe. \textbf{43} (1983), 23--33.



\end{thebibliography}
\end{document}